\documentclass[12pt]{article}
\usepackage{amssymb, graphicx}

\newtheorem{lemma}{Lemma}[section]

\newtheorem{constr}{Construction}[section]

\newcommand{\nl}{\par \bigskip \noindent}

\title{Transitive decompositions of graphs and their links with geometry and origami}
\author{Geoffrey Pearce}
\begin{document}
\maketitle
\begin{abstract}
A transitive decomposition of a graph is a partition of the edge or arc set giving a set of subgraphs which are preserved and permuted transitively by a group of automorphisms of the graph.  In this paper we give some background to the study of transitive decompositions and highlight a connection with partial linear spaces.  We then describe a simple method for constructing transitive decompositions using graph quotients, and we show how this may be used in an application to modular origami.
\end{abstract}
\section{Introduction}
The idea of a transitive decomposition of a graph has been around for a long time (since the 1970s), and has appeared in a number of different guises; but only recently has it been studied extensively in its own right.  Any graph can be broken down into a set of subgraphs by partitioning the edge set, and this is known as a {\em decomposition} of the graph.  This generalises the more widely-known idea of a graph factorisation, which has the extra requirement that the subgraphs are {\em spanning}; that is, each vertex of the whole graph is incident with some edge in each of the subgraphs.  

During the 1970s and 1980s, considerable energy was devoted to the study of graph factorisations in which the subgraphs are pairwise isomorphic (see \cite{HarRob}).  A transitive decomposition can be thought of as a `special case of a generalisation' of this isomorphic factorisation idea -- it is a generalisation because it allows subgraphs whose vertices do not span the whole graph, and a special case because (as we explain in a moment) it requires that the subgraphs are isomorphic `with respect to the symmetry of the graph'.  

Formally, a transitive decomposition is a pair $(\Gamma,{\mathcal P})$ where $\Gamma$ is a graph (with vertex set $V\Gamma$ and edge set $E\Gamma$), ${\mathcal P}$ is a partition of $E\Gamma$, such that there is a group $G$ of automorphisms of $\Gamma$ (in other words, a group of permutations of $V\Gamma$ which leaves the edge set invariant) satisfying the following two conditions:
\begin{itemize}
\item[(i)] given any part $P$ in ${\mathcal P}$ and any automorphism in $G$, $P$ is mapped by the automorphism either wholly to itself or wholly to a different part in ${\mathcal P}$; and
\item[(ii)] given any two parts in the partition, there exists an automorphism in $G$ mapping the first wholly to the second.
\end{itemize}
When this holds we call $(\Gamma, {\mathcal P})$ a $G$-transitive decomposition.  Each part $P$ in the partition corresponds to a subgraph $\Gamma_P$ whose edge set is $P$ and whose vertex set consists of those vertices of $\Gamma$ incident with edges in $P$.  Thus the first of these conditions is equivalent to requiring that the group $G$ preserves the set of subgraphs in the decomposition; and the second is equivalent to requiring that $G$ is transitive on the set of subgraphs.  (In general a group $G$ is {\em transitive} on a set $\Omega$ if for all $\omega, \omega' \in \Omega$, there is an element of $G$ mapping $\omega$ to $\omega'$.)  Note that this second condition forces the subgraphs to be pairwise isomorphic, since an automorphism mapping one subgraph to another induces an isomorphism between the two.  A transitive decomposition may also be defined more generally with ${\mathcal P}$ a partition of the arcs (ordered pairs of adjacent vertices) of $\Gamma$.  However for simplicity we will focus on the case where ${\mathcal P}$ is an edge partition.

Within algebraic graph theory, special classes of transitive decompositions have been studied in a number of isolated cases.  Most of this research interest has focussed on classes of isomorphic factorisations with an automorphism group acting transitively on the subgraphs (see for example \cite{burfra, camkor, PraLi2}).  However, some researchers have shown interest in the more general notion of a transitive  decomposition.  For example, in \cite{Sibley} Thomas Sibley gives a characterisation of `edge-coloured graphs with a two-transitive automorphism group'.  This is an alternative way of describing decompositions of complete graphs whose automorphism group acts two-transitively on vertices (a graph is {\em complete} if every pair of vertices is an edge).

One of the main reasons for the current interest in transitive decompositions is the large number of connections they have with well-known structures in combinatorics and geometry.  We explain some of these connections in the next section.  In Section \ref{constructions} we describe a method for constructing transitive decompositions using graph quotients, and in Section \ref{origami} we outline an interesting application of transitive decompositions to modular origami.

\section{Connections with other structures} \label{connections} 
Transitive decompositions have some of their most interesting connections with structures that traditionally lie outside of graph theory.  Many of these connections do not occur in the context of graph {\em factorisations} because they require that the vertex sets of the subgraphs are proper subsets of the vertex set of the whole graph.  As an example, we describe below the equivalence of line-transitive partial linear spaces with certain transitive decompositions.

A {\em partial linear space} is a set of points together with a set of (at least two) lines. Each line is a subset of points, and every pair of points lies in at most one line.  We will denote the point set by ${\mathcal V}$, the line set by ${\mathcal L}$, and the partial linear space itself by the pair $({\mathcal V},{\mathcal L})$.  A partial linear space is said to be {\em line transitive} if there is a group of permutations of the points which preserves and transitively permutes the lines.  Lemma \ref{ltpls} shows that every line transitive partial linear space gives rise to a transitive decomposition of a certain type, and conversely that every transitive decomposition of this type gives rise to a line transitive partial linear space.  Given a permutation $g$ and an object $x$ permuted by $g$ we write $x^g$ to denote the image of $x$ under $g$.

\begin{lemma}\label{ltpls}\ \vspace{1pt}
\begin{itemize}
\item[(i)] Let $({\mathcal V},{\mathcal L})$ be a line transitive partial linear space, and suppose that $G$ is a group of permutations of ${\mathcal V}$ which preserves and acts transitively on ${\mathcal L}$.  Let $\Gamma$ be the graph with vertex set ${\mathcal V}$ and edges $\{\alpha,\beta\}$ whenever there exists $\ell \in {\mathcal L}$ with $\alpha$ and $\beta$ both in $\ell$.  For each $\ell \in {\mathcal L}$, let $P_{\ell}$ be the set of all unordered pairs of distinct elements of ${\ell}$, and let ${\mathcal P} = \{P_{\ell} \, | \, {\ell} \in {\mathcal L}\}$.  Then $(\Gamma,{\mathcal P})$ is a $G$-transitive decomposition, and each $\Gamma_{P_{\ell}}$ is a complete subgraph of $\Gamma$.
\item[(ii)] Let $(\Gamma,{\mathcal P})$ be a $G$-transitive decomposition such that for each $P \in {\mathcal P}$, the subgraph $\Gamma_P$ is a complete subgraph of $\Gamma$.  Let ${\mathcal V} = V\Gamma$, and let ${\mathcal L} = \{V\Gamma_P \, | \, P \in {\mathcal P}\}$.  Then $({\mathcal V},{\mathcal L})$ is a line transitive partial linear space.
\end{itemize}
\end{lemma}
{\em Proof.} To prove part (i), note first that since each pair of points in ${\mathcal V}$ lies in at most one line in ${\mathcal L}$, each edge of $\Gamma$ lies in at most one part of ${\mathcal P}$.  The definition of $E\Gamma$ implies that $\bigcup_{\ell \in {\mathcal L}} P_\ell$ is the whole of $E\Gamma$, and so ${\mathcal P}$ defines a partition of $E\Gamma$.  To see that ${\mathcal P}$ is $G$-invariant, observe that if an element $g \in G$ maps an edge in a part $P_{\ell}$ to an edge in a part $P_{{\ell}'}$, then ${\ell}^g$ contains two points of ${\ell}'$.  Hence ${\ell}^g = {\ell}'$, since each pair of points lies in at most one line, and it follows that $P_{\ell}^g = P_{{\ell}'} \in {\mathcal P}$.  Furthermore, $G$ acts transitively on ${\mathcal P}$, since $G$ acts transitively on lines.  Since each part $P_{\ell}$ is the set of all pairs of distinct elements of ${\ell}$, it is clear that $\Gamma_{P_{\ell}}$ is a complete subgraph of $\Gamma$.

To prove part (ii) we note that since every edge $\{\alpha,\beta\}$ of $\Gamma$ lies in exactly one part in ${\mathcal P}$, every pair of points of ${\mathcal V}$ lies in at most one line in ${\mathcal L}$.  For any $g$ in $G$ we have $(V\Gamma_P)^g = V\Gamma_{P^g}$ and it follows that $G$ preserves ${\mathcal L}$.  Furthermore, $G$ acts transitively on ${\mathcal P}$ and hence on ${\mathcal L}$, and so $({\mathcal V},{\mathcal L})$ is a line transitive partial linear space.

\nl
% \begin{figure} [h]
% \begin{center}
% \begin{tabular}{ccc}
% \includegraphics[height=1.4in]{bwlinetranspls.eps} &
% \includegraphics[height=1.4in]{arrow.eps} &
% \includegraphics[height=1.4in]{bwtransdecomp.eps}
% \end{tabular}
% \end{center}
% \caption{A line transitive partial linear space (left), and the corresponding graph (right).  The partial linear space has 3 lines: a solid line, a dotted line and a zigzag line, and these correspond respectively to the solid, dotted and zigzag subgraphs on the right, giving a transitive decomposition of the graph.  Note that each such subgraph in the decomposition is a complete graph with 3 vertices.} \label{ltpls}
% \end{figure}

Several other geometrical structures exhibit similar connections with certain families of transitive decompositions.  These include linear spaces, cyclic Hamiltonian cycle systems \cite{burfra}, symmetric association schemes \cite{bannai}, and symmetric graph designs \cite{cam}.

\section{Using quotients to construct transitive decompositions} \label{constructions}
As is the case for many structures in mathematics, it is often possible to construct transitive decompositions by `multiplying' or `expanding' known examples.  In \cite{Pearce}, the notion of taking a `product' of transitive decompositions is explored in detail.  The constructions described therein have important applications in the study of transitive decompositions in which the group $G$ is a {\em primitive group of grid type} with {\em rank $3$} (this means that $G$ acts transitively on the vertices, on the arcs, and on the `non-arcs' of $\Gamma$).

Here we consider another useful method for constructing transitive decompositions which involves the graph-theoretical concept of a quotient.  Let $\Gamma$ be a graph and suppose that $G$ is a group of automorphisms of $\Gamma$ which acts transitively on $V\Gamma$.  Suppose also that ${\mathcal B}$ is a partition of $V\Gamma$ which is $G$-invariant (this means that ${\mathcal B}$ as a set of subsets is left unchanged by each permutation in $G$).  The group $G$ therefore permutes the elements of ${\mathcal B}$, and we write $G^{\mathcal B}$ for the group of permutations of ${\mathcal B}$ induced by this action of $G$.

Let $\Gamma_{\mathcal B}$ be the graph with vertex set ${\mathcal B}$ and edges $\{B_1,B_2\}$ whenever there is an edge of $\Gamma$ between $B_1$ and $B_2$; that is, whenever there exists $\alpha \in B_1$ and $\beta \in B_2$ with $\{\alpha,\beta\} \in E\Gamma$.  Then $\Gamma_{\mathcal B}$ is called the {\em imprimitive quotient} of $\Gamma$ with respect to ${\mathcal B}$.  We show now that it is sometimes possible to construct a transitive decomposition of $\Gamma$ given a transitive decomposition of $\Gamma_{\mathcal B}$.

\begin{constr} \label{cover} {\em
Let $\Gamma$ be a graph and let ${\mathcal B}$ be a $G$-invariant partition of $V\Gamma$ such that each $B_i \in {\mathcal B}$ consists of pairwise non-adjacent vertices of $\Gamma$.  Suppose that $(\Gamma_{\mathcal B},{\mathcal Q})$ is a $G^{\mathcal B}$-transitive decomposition.  For each $Q$ in ${\mathcal Q}$ define a subset $P_Q$ of $E\Gamma$ by
$$P_Q = \{ \{\alpha,\beta\} \, | \, \{B_i,B_j\} \in Q \; \mbox{with} \; \alpha \in B_i \; \mbox{and} \; \beta \in B_j \; \mbox{and} \; \{\alpha,\beta\} \in E\Gamma\}$$
and let ${\mathcal P} = \{P_Q \, | \, Q \in {\mathcal Q}\}$.
}
\end{constr}

\begin{lemma} \label{coverlem}
$(\Gamma,{\mathcal P})$ is a $G$-transitive decomposition.
\end{lemma}
{\em Proof.} First, observe that no pair $\{B_i,B_j\}$ occurs in more than one part in ${\mathcal Q}$, and so no edge of $\Gamma$ can appear in more than one part of ${\mathcal P}$.  Also since ${\mathcal Q}$ is a partition of $E\Gamma_{\mathcal B}$, and since no vertices in any $B_i$ are adjacent in $\Gamma$, it follows from the definition of the parts in ${\mathcal P}$ that each edge of $\Gamma$ lies in some part in ${\mathcal P}$.  Thus ${\mathcal P}$ is a partition of $E\Gamma$.

Let $P_Q \in {\mathcal P}$.  For any $g \in G$, the image of $P_Q$ under $g$ consists of all edges which lie between $B_i^g$ and $B_j^g$ for all $\{B_i,B_j\} \in Q$. This means that $(P_Q)^g = P_{Q^g}$; and since ${\mathcal Q}$ is left invariant by $G^{\mathcal B}$, it follows that ${\mathcal P}$ is left invariant by $G$.  Moreover, since $G^{\mathcal B}$ is transitive on ${\mathcal Q}$, $G$ is transitive on ${\mathcal P}$.

\nl
Similar constructions to this one have proved useful in characterising certain transitive decompositions with $G$ an imprimitive rank $3$ group and $\Gamma$ a complete multipartite graph.
\section{An application to modular origami} \label{origami}
The traditional Japanese art of paper-folding known as origami has become a very popular recreational activity throughout the world, appealing to both children and adults.  Over the last century, a number of variations and off-shoots of the traditional artform have appeared, and one of the most interesting of these is known as `modular' origami (or sometimes `unit' origami).  This usually involves building large geometric structures from several smaller, individually folded `modules' which are structurally identical and fitted together without adhesive (see \cite{Simon} for some detailed examples and instructions).  Such structures can be built with as few as 6 modules, or as many as 900 (or even 1720 in a model built by the author and a friend during their idle undergraduate years).  Models created in this way are often spectacular and decorative, especially when different colours are used for the modules.  Trying to find the best ways of using different colours gives rise to some interesting questions not only in aesthetics, but also in mathematics.  There are a couple of colouring conditions which usually produce spectacular models:
\begin{itemize}
\item[(i)] no two modules of the same colour are joined;
\item[(ii)] the distribution of colours over the model is `symmetric' or `regular' in some sense.
\end{itemize}
It turns out that most modular origami models have an underlying graph structure: a module can be thought of as an edge, and the point at which two or more modules meet can be thought of as a vertex.  These underlying graphs often correspond to geometric solids, with the graph of the dodecahedron arising especially frequently.  This graph theoretic model gives an alternative way of understanding the above conditions: a colouring satisfying condition (i) corresponds to an edge-decomposition of the underlying graph in which each subgraph has valency $1$ (in other words, no two edges in a subgraph share a vertex), and the problem of satisfying condition (ii) can be solved by finding a transitive decomposition of the underlying graph.  Solving both together therefore involves finding a transitive decomposition in which each subgraph has valency $1$; or in other words, a transitive $1$-decomposition.

\subsection{A modular origami colouring derived from a transitive decomposition of the dodecahedron}
A dodecahedron has two types of symmetries: those which can be achieved by rotation in space, and those which involve inverting the object through itself.  Let $\Gamma$ be the graph of the dodecahdron, and let $G$ denote the group of automorphisms of $\Gamma$ corresponding to rotation symmetries of the dodecahedron.  It is not hard to see that any vertex of a dodecahedron may be moved to any other by some rotation through space, and this means that $G$ acts transitively on $V\Gamma$.  Also, it can be seen intuivitely that any rotational symmetry moves any antipodal pair of vertices wholly to another (or the same) antipodal pair.  This means that if we take ${\mathcal B}$ to be the set of all antipodal pairs of $\Gamma$, then the group $G$ leaves ${\mathcal B}$ invariant as a partition of $V\Gamma$.  Thus we may construct the imprimitive quotient $\Gamma_{\mathcal B}$.

\begin{figure} [h]
\begin{center}
\includegraphics[height=3.0in]{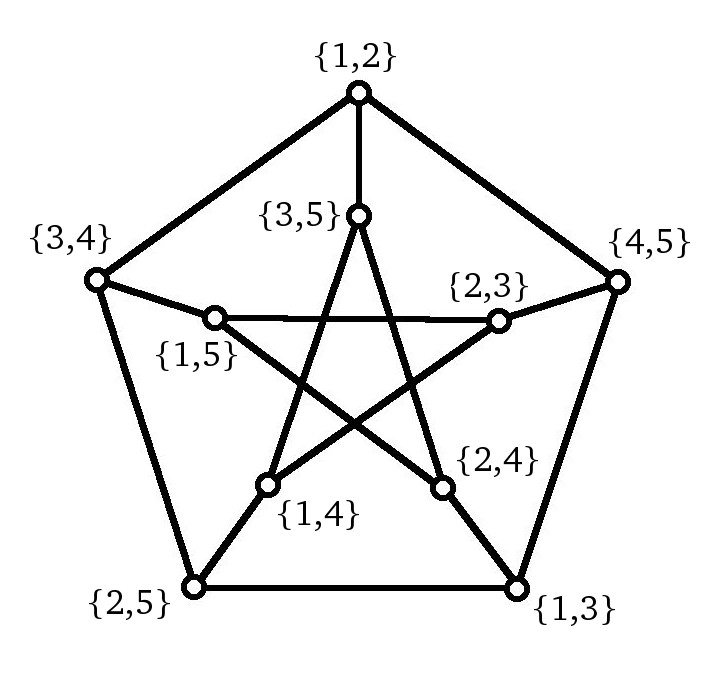} 
\end{center}
\caption{The Petersen graph with vertices labelled as described above.} \label{Petersen}
\end{figure}

It is well-known that $\Gamma_{\mathcal B}$ is isomorphic to the Petersen graph; and furthermore that the Petersen graph admits a labelling of the vertices by unordered pairs of elements from the set $\{1,2,3,4,5\}$ such that $\{a,b\}$ is adjacent to $\{c,d\}$ whenever $\{a,b\} \cap \{c,d\} = \varnothing$ (see Figure \ref{Petersen}). 

This labelling gives us an easy way to construct a transitive decomposition of the Petersen graph.  We will use this together with Construction \ref{cover} to find a transitive $1$-decomposition of the dodecahedron, and thereby find a modular origami colouring satisfying the colouring conditions (i) and (ii).

Let $\Delta$ denote the Petersen graph with vertices labelled as described above.  It is well-known that $A_5$, the alternating group of degree $5$, acts as a group of automorphisms of $\Delta$ by permuting the vertices according to $\{a,b\}^g = \{a^g,b^g\}$ for all $\{a,b\} \in V\Delta$ and $g \in A_5$.  Indeed we may identify the action induced by $G$ on ${\mathcal B}$ with this action of $A_5$ on $V\Delta$.

For each $a \in \{1,\ldots,5\}$, define $Q_a$ to be the subset of $E\Delta$ consisting of the three edges $\{\{b,c\},\{d,e\}\}$ such that $a \not\in \{b,c,d,e\}$, and define
$${\mathcal Q} := \{Q_a \, | \, a \in \{1,\ldots,5\}\}.$$
Clearly ${\mathcal Q}$ is a partition of $E\Delta$ containing five parts.

\begin{lemma}
$(\Delta,{\mathcal Q})$ is an $A_5$-transitive $1$-decomposition.
\end{lemma}
{\em Proof.}
First observe that any edge $\{\{b,c\},\{d,e\}\}^g$ in $(Q_a)^g$ is such that $a^g \not\in \{b^g,c^g,d^g,e^g\}$, and so $(Q_a)^g = Q_{a^g}$, which shows that $G$ preserves the partition ${\mathcal Q}$.  Given that $A_5$ acts in this way on ${\mathcal Q}$, it follows from the transitivity of $A_5$ on $\{1,\ldots,5\}$ that $A_5$ is also transitive on ${\mathcal Q}$.  Hence $(\Delta,{\mathcal Q})$ is an $A_5$-transitive decomposition.

Now suppose that two different edges $\{\{b,c\},\{d,e\}\}$ and $\{\{b',c'\},\{d',e'\}\}$ in $Q_a$ share a vertex.  Without loss of generality we may suppose that $\{b,c\} = \{b',c'\}$.  Then $\{d,e\} \neq \{d',e'\}$, and since also $\{b,c\} \cap \{d,e\} = \varnothing$ and $\{b,c\} \cap \{d',e'\} = \varnothing$, it follows that $\{b,c,d,e,d',e'\}$ contains the five distinct elements of $\{1,\ldots,5\}$.  Hence $a \in \{b,c,d,e,d',e'\}$; but this contradicts the definition of $Q_a$, and so $\{\{b,c\},\{d,e\}\}$ and $\{\{b',c'\},\{d',e'\}\}$ cannot share a vertex.  Hence each subgraph in the decomposition has valency $1$.

\nl
Now, since two antipodal vertices of $\Gamma$ can never be adjacent, it follows that each subset in ${\mathcal B}$ consists of pairwise non-adjacent vertices.  Hence, identifying $\Delta$ with $\Gamma_{\mathcal B}$ and $A_5$ with $G^{\mathcal B}$, we may now apply Construction \ref{cover} and Lemma \ref{coverlem} to obtain a $G$-transitive decomposition $(\Gamma,{\mathcal P})$.  If some part $P \in {\mathcal P}$ contained two edges sharing a vertex, then there would be a part $Q \in {\mathcal Q}$ containing two edges also sharing a vertex, and this is not the case; hence $(\Gamma,{\mathcal P})$ is a $G$-transitive $1$-decomposition.

\nl
If we now assign a different colour to each of the parts $P_a$ and construct the corresponding modular origami model according to this colour scheme, we find that no two modules of the same colour are joined, and that the distribution of colours over the model is symmetrical.  This results in a very decorative colouring with five different colours.
\nl
\begin{figure} [h]
\begin{center}
\includegraphics[height=2.5in]{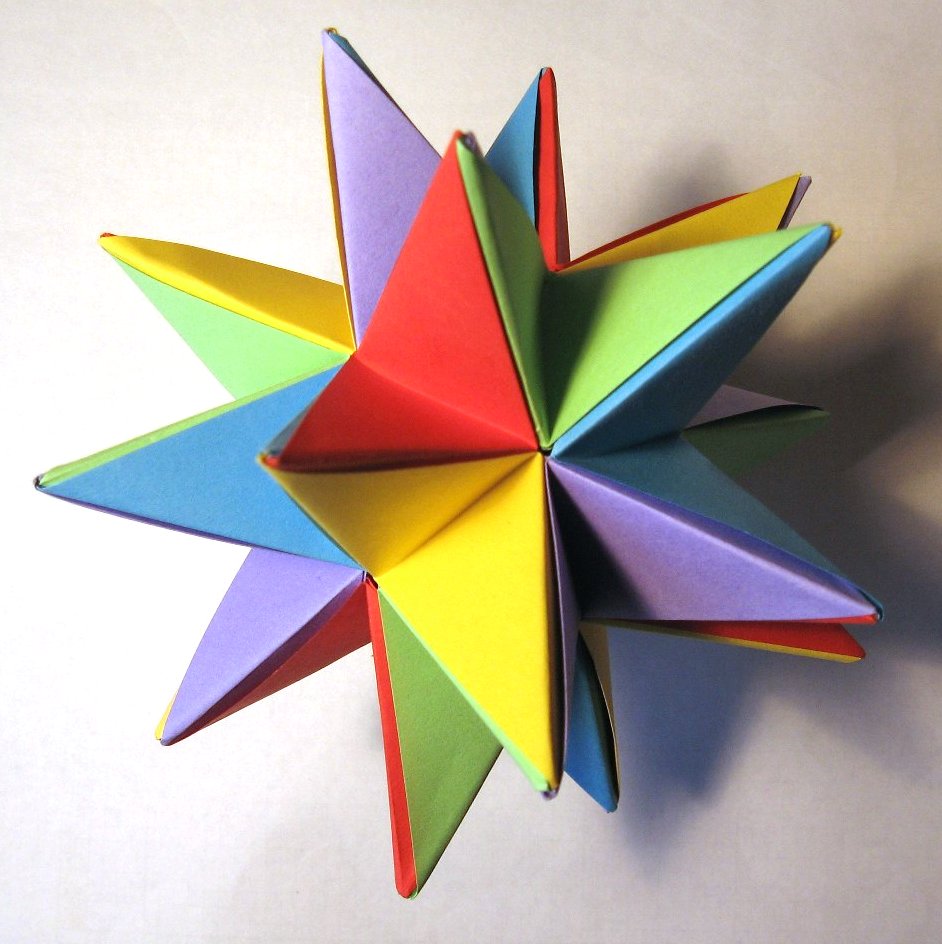} 
\end{center}
\caption{A model based on the transitive decomposition described above.} \label{model}
\end{figure}
\section{Acknowledgements}
The author gratefully acknowledges assistance from his Ph.D. supervisors, Cheryl E. Praeger and John Bamberg, in the writing of this paper.  The author was supported by an Australian Postgraduate Award.
\bibliography{gazette}
\bibliographystyle{plain}
\end{document}